\documentclass{elsart}
\usepackage{amsfonts}
\usepackage{amsmath}

\input{tcilatex}
\begin{document}

\begin{frontmatter}%

\title{Matrices about the lattice paths}%

\author{Jishe Feng}%
\footnote{Corresponding author. E-mail: gsfjs6567@126.com}%

\collab{}%

\address{Department of Mathematics, Longdong University,  Qingyang,  Gansu,  745000,  China}%

\begin{abstract}
From the matrix point of view, we use the recursion to discuss four
combinatorial numbers in terms of the integer lattice paths, this is
different from Andr\'{a}'s method \cite{Andra}. We give four tables and
matrices, and their relations, and give their matrix decompositions. By
their relations, we give the explicit formulas of the special combinatorial
numbers.

\textit{AMS} classification: 05A15

\end{abstract}%

\begin{keyword}
Lattice path; combinatorial number; recurrence formula; matrix decomposition%
\end{keyword}%

\end{frontmatter}%

\section{Introduction}

Many special integers can be interpreted to the numbers of path on the
integer lattice in the coordinate plane, such as Catalan number and Schr\H{o}%
der number can be interpreted to the number of lattice paths that start at $%
(0,0)$, end at $(n,n)$, contain no points above the line $y=x$, and are
composed only of horizontal step $(0,1)$, vertical step $(1,0)$, and up
diagonal step $(1,1)$, i.e.,$\rightarrow $ ,$\uparrow $ , and $\nearrow $ 
\cite{Richard}. In paper [2], the authors discuss the ballot problem,
Catalan number and their generalization in terms of $p$-good paths on the
integer lattice, in order to get the formula, they used Andr\'{a}'s
reflection method [1,2] to subtract the numbers of the bad path from the
total number to obtain the number of good path. Renault [3,4] described the
first combinatorial proof as given by Andr\'{a} \cite{Andra} and claimed
that the reflection method typically misattributed to Andr\'{a}. Goulden and
Serrano \cite{Goulden} provide a direct geometric bijection for the number
of lattice paths that never go below the line $y=kx$ for a positive integer $%
k$, it uses rotation instead of reflect.

From the historical point of view, the oldest are known as Pascal matrix $P$%
, and its $n$ order minor $P_{n},$ they are as follows.%
\begin{equation*}
P=\left( 
\begin{array}{ccccccc}
\binom{0}{0} &  &  &  &  &  &  \\ 
\binom{1}{0} & \binom{1}{1} &  &  &  &  &  \\ 
\binom{2}{0} & \binom{2}{1} & \binom{2}{2} &  &  &  &  \\ 
\vdots & \vdots & \vdots & \vdots &  &  &  \\ 
\binom{n-1}{0} & \binom{n-1}{1} & \binom{n-1}{2} & \cdots & \binom{n-1}{n-1}
&  &  \\ 
\vdots & \vdots & \vdots & \vdots & \vdots &  & \ddots%
\end{array}%
\right) ,P_{n}=\left( 
\begin{array}{ccccc}
\binom{0}{0} &  &  &  &  \\ 
\binom{1}{0} & \binom{1}{1} &  &  &  \\ 
\binom{2}{0} & \binom{2}{1} & \binom{2}{2} &  &  \\ 
\vdots & \vdots & \vdots & \vdots &  \\ 
\binom{n-1}{0} & \binom{n-1}{1} & \binom{n-1}{2} & \cdots & \binom{n-1}{n-1}%
\end{array}%
\right)
\end{equation*}
Various types of Pascal matrices are investigated in [7-10]. Shapiro \cite%
{Shapiro} introduced a number triangle with its first column entries
containing the Catalan numbers, Barcucci \cite{Barcucci} set up the Catalan
triangle, which is a number triangle with entries equal to $\frac{n-k+1}{n+1}%
\binom{n+k}{k}.$

In this paper, from the matrix point of view, we use the recursion to
discuss four combinatorial numbers in terms of the integer lattice paths,
this is different from Andr\'{a}'s method \cite{Andra}. It is the purpose of
this paper to discuss four combinatorial numbers about the lattice path, and
give four tables and matrices which are different from the Toeplitz matrices
which possess numbers of various type arranged on the main diagonal and
below \cite{Stefan}, and give their matrix decompositions. By their
relations, we give the explicit formulas of the special combinatorial
numbers.

\section{Rectangular lattice path and combinatorial number}

For convenience to express as a matrix, we choose the vertical down line as $%
x$-axis and the horizontal line as $y$-axis, so the vertical step $\uparrow $
changes to $\downarrow $, up diagonal step $\nearrow $ changes to down
diagonal step $\searrow $. For any positive integer $k,r$, A rectangular
lattice path \cite{Richard} from $(0,0)$ to $(k,r)$ is a path from $(0,0)$
to $(k,r)$ that is made up of horizontal steps $\rightarrow $ and vertical
steps $\downarrow $. Let $b(k,r)$ denote the number of the rectangular
lattice path from $(0,0)$ to $(k,r)$.

In order to get the explicit expression of the number of the rectangular
lattice path, we can analyze it step by step. For each point on $x$-axis or $%
y$-axis, there is only one path from $(0,0)$ to them, it is only $%
\rightarrow $ or $\downarrow $, namely $b(k,0)=b(0,r)=1$. For the point $%
(1,1)$, that is either from the left point $(1,0)$ or from the above point $%
(0,1)$ to it, thus the number $b(1,1)$ is the sum of $b(1,0)$ and $\ b(0,1)$%
, that is $b(1,1)=b(1,0)+b(0,1)$. In general for the point $(k,r)$, the last
step of the rectangular lattice path from $(0,0)$ to $(k,r)$ is either from $%
(k,r-1)$ by $\rightarrow $ or from $(k-1,r)$ by $\downarrow $ to it, so
there is 
\begin{equation}
b(k,r)=b(k,r-1)+b(k-1,r).  \label{abc1}
\end{equation}

By tabulating the number of rectangular lattice paths from $(0,0)$ to each
point$,$ we obtain a symmetric number table (see table 1) and symmetric
matrix $Q.$

$\ \ 
\begin{tabular}{c|cccccccc}
\multicolumn{9}{l}{Table 1} \\ \hline
$k\backslash r$ & $0$ & $1$ & $2$ & $3$ & $4$ & $5$ & $6$ & $\cdots $ \\ 
\hline
$0$ & $1$ & $1$ & $1$ & $1$ & $1$ & $1$ & $1$ & $\cdots $ \\ 
$1$ & $1$ & $2$ & $3$ & $4$ & $5$ & $6$ & $7$ & $\cdots $ \\ \cline{6-6}
$2$ & $1$ & $3$ & $6$ & $10$ & \multicolumn{1}{|c}{$15$} & 
\multicolumn{1}{|c}{$21$} & $28$ & $\cdots $ \\ \cline{5-5}
$3$ & $1$ & $4$ & $10$ & \multicolumn{1}{|c}{$20$} & $35$ & 
\multicolumn{1}{|c}{$56$} & $84$ & $\cdots $ \\ \cline{5-6}
$4$ & $1$ & $5$ & $15$ & $35$ & $70$ & $126$ & $210$ & $\cdots $ \\ 
$5$ & $1$ & $6$ & $21$ & $56$ & $126$ & $252$ & $462$ & $\cdots $ \\ 
$6$ & $1$ & $7$ & $28$ & $84$ & $210$ & $462$ & $924$ & $\cdots $ \\ 
$\vdots $ & $\vdots $ & $\vdots $ & $\vdots $ & $\vdots $ & $\vdots $ & $%
\vdots $ & $\vdots $ & 
\end{tabular}%
$,

$Q=\left[ 
\begin{tabular}{cccccccc}
$1$ & $1$ & $1$ & $1$ & $1$ & $1$ & $1$ & $\cdots $ \\ 
$1$ & $2$ & $3$ & $4$ & $5$ & $6$ & $7$ & $\cdots $ \\ 
$1$ & $3$ & $6$ & $10$ & $15$ & $21$ & $28$ & $\cdots $ \\ 
$1$ & $4$ & $10$ & $20$ & $35$ & $56$ & $84$ & $\cdots $ \\ 
$1$ & $5$ & $15$ & $35$ & $70$ & $126$ & $210$ & $\cdots $ \\ 
$1$ & $6$ & $21$ & $56$ & $126$ & $252$ & $462$ & $\cdots $ \\ 
$1$ & $7$ & $28$ & $84$ & $210$ & $462$ & $924$ & $\cdots $ \\ 
$\vdots $ & $\vdots $ & $\vdots $ & $\vdots $ & $\vdots $ & $\vdots $ & $%
\vdots $ & 
\end{tabular}%
\right] $

\begin{thm}
The number of rectangular lattice paths from $(0,0)$ to $(k,r)$ equals the
binomial coefficient, namely $b(k,r)=$ $\binom{k+r}{r}.$

\begin{pf}
A rectangular lattice path from $(0,0)$ to $(k,r)$ is uniquely determined by
its sequence of $k$ horizontal steps $\rightarrow $ and $r$ vertical steps $%
\downarrow $, and every such sequence determines a rectangular lattice path
from $(0,0)$ to $(k,r)$. Hence, the number of paths equals the number of
permutations of $k+r$ objects of which $k$ are $\rightarrow $ and $s$ are $%
\downarrow $, that is the result \cite{Richard}.
\end{pf}
\end{thm}

According to Pascal's formula $\binom{n-1}{m-1}+\binom{n-1}{m}=\binom{n}{m}$%
, and the formula ($\ref{abc1}),$ we can get the propositions as follows.

\begin{prop}
Let $Q(p,q)$ be the element in $p$th row and $q$th column of the matrix $Q$,
then every element of the symmetric matrix $Q$ in first row and first column
is $1$, and every other element is the sum of its above element and left
element. 
\begin{equation}
Q(p,q)=Q(p,q-1)+Q(p-1,q),  \label{b1}
\end{equation}%
\begin{equation}
Q(p,q)=\underset{i=0}{\overset{q}{\sum }}Q(p-1,i).  \label{b2}
\end{equation}
\end{prop}

\begin{prop}
$Q_{n}$ can be decomposed as the product of Pascal matrix and its
transposition. 
\begin{equation}
Q_{n}=P_{n}P_{n}^{T}  \label{b3}
\end{equation}
\end{prop}

This property has been in paper \cite{Brawer}. By carrying out the
multiplication of the matrices $P_{n}$ and $P_{n}^{T}$, we can obtain the
two identities

\begin{prop}
\begin{equation}
\underset{k=0}{\overset{\min (i,j)-1}{\sum }}\binom{i-1}{k}\binom{j-1}{k}=%
\binom{i+j-2}{i-1};  \label{a1}
\end{equation}%
\begin{equation}
\underset{k=0}{\overset{i-1}{\sum }}\binom{i-1}{k}^{2}=\binom{2i-2}{i-1},
\label{a2}
\end{equation}%
where $i=1,2,3,\cdots .$
\end{prop}

\section{Subdiagonal rectangular lattice path and Catalan number}

We now consider rectangular lattice paths from $(0,0)$ to $(p,q)$ that are
restricted to lie on or below the line $y=x$ in the coordinate plane
whenever $p\geq q$. Brualdi \cite{Richard} call such paths subdiagonal
rectangular lattice paths. Let $C(p,q)$ denote the number of the subdiagonal
rectangular lattice path from $(0,0)$ to $(p,q)$.

We can also analyze it step by step. For each point on $x$-axis, there is
only one path from $(0,0)$ to them, it only $\downarrow $, namely $C(p,0)=1$%
. For the point $(1,1)$, it is only from the left point $(1,0)$ to it, thus
the number $C(1,1)$ is $C(1,0)$, that is $C(1,1)=C(1,0)$. For the point $%
(2,1)$, its last step is either from the left point $(2,0)$ or from the
above point $(1,1)$ to it, namely, $C(2,1)=C(2,0)+C(1,1)$, and so on, we can
obtain a triangle, which is same as the triangle in \cite{Barcucci} (see
table 2 ).

$%
\begin{tabular}{c|ccccccccc}
\multicolumn{10}{l}{Table 2} \\ \hline
$p\backslash q$ & $0$ & $1$ & $2$ & $3$ & $4$ & $5$ & $6$ & $7$ & $\cdots $
\\ \hline
$0$ & $1$ &  &  &  &  &  &  &  &  \\ 
$1$ & $1$ & $1$ &  &  &  &  &  &  &  \\ 
$2$ & $1$ & $2$ & $2$ &  &  &  &  &  &  \\ 
$3$ & $1$ & $3$ & $5$ & $5$ &  &  &  &  &  \\ \cline{6-6}
$4$ & $1$ & $4$ & $9$ & $14$ & \multicolumn{1}{|c}{$14$} & 
\multicolumn{1}{|c}{} &  &  &  \\ \cline{5-5}
$5$ & $1$ & $5$ & $14$ & \multicolumn{1}{|c}{$28$} & $42$ & 
\multicolumn{1}{|c}{$42$} &  &  &  \\ \cline{5-6}
$6$ & $1$ & $6$ & $20$ & $48$ & $90$ & $132$ & $132$ &  &  \\ 
$7$ & $1$ & $7$ & $27$ & $75$ & $165$ & $297$ & $429$ & $429$ &  \\ 
$\vdots $ & $\vdots $ & $\vdots $ & $\vdots $ & $\vdots $ & $\vdots $ & $%
\vdots $ & $\vdots $ & $\vdots $ & $\cdots $%
\end{tabular}%
$

Because the restriction that it must lie on or below the line $y=x$, the
element in diagonal $C(n,n)$ equals the element $C(n,n-1)$. Compare the
lower triangle of Table 1 and Table 2, we find that the element in Table 2 $%
C(p,q)$ equal the element in Table 1 $b(p,q)$ subtract the element $%
b(p+1,q-1)$, that is $C(p,q)=b(p,q)-b(p+1,q-1)$. This is derived from the
symmetric of the Table 1, $b(p+1,q-1)$ is the number of subdiagonal
rectangular paths above the line $y=x$, it is sure that it must be
subtracted. By Theorem 1, $b(p,q)=\binom{p+q}{p}$, we can get%
\begin{eqnarray*}
C(p,q) &=&b(p,q)-b(p+1,q-1)=\binom{p+q}{q}-\binom{p+1+q-1}{q-1} \\
&=&\binom{p+q}{q}-\binom{p+q}{q-1}=\frac{p-q+1}{p+1}\binom{p+q}{q}.
\end{eqnarray*}%
Thus, we can get the following theorem.

\begin{thm}
For $C(p,q)$, The number of subdiagonal rectangular lattice paths from $%
(0,0) $ to $(p,q)$, there is \cite{Richard}%
\begin{equation}
C(p,q)=\frac{p-q+1}{p+1}\binom{p+q}{q}.  \label{r3}
\end{equation}
\end{thm}

There are the same formula as (\ref{b1}) and (\ref{b2}) 
\begin{equation}
C(p,q)=C(p-1,q)+C(p,q-1),  \label{s3}
\end{equation}%
\begin{equation}
C(p+1,q)=\underset{k=0}{\overset{q}{\sum }}C(p,k).  \label{s4}
\end{equation}

It is easy found that diagonal elements and its left elements are \textit{%
Catalan numbers}, this means that they are the number of subdiagonal
rectangular lattice paths from $(0,0)$ to $(n,n)$. From the formula (\ref{r3}%
), one obtains $C(n,n)=\frac{1}{n+1}\binom{2n}{n}$, this is the well-known
Catalan number $C_{n}$. By the formula (\ref{s3}) and (\ref{s4}), we can get
two identities and the recurrence formula about Catalan number as follows:%
\begin{equation}
C_{n}=C(n,n-1)=\frac{1}{n}\binom{2n-2}{n-1}+\frac{3}{n+1}\binom{2n-2}{n-2},
\label{s5}
\end{equation}%
\begin{equation}
C_{n}=C(n,n)=\underset{k=0}{\overset{n-1}{\sum }}\frac{n-k}{n}\binom{n+k-1}{k%
}=\frac{1}{n+1}\binom{2n}{n},  \label{s6}
\end{equation}%
\begin{equation}
C_{n}=C(n,n)=\frac{2}{n+1}\binom{2n-1}{n-1}=\frac{2(2n-1)}{n+1}C_{n-1}.
\label{s66}
\end{equation}

\section{HVD-lattice path and Delanny number}

A $HVD$-lattice path \cite{Richard} is a rectangular lattice path from $%
(0,0) $ to $(p,q)$ is a path from $(0,0)$ to $(p,q)$ that is made up of
horizontal steps $\rightarrow $, vertical steps $\downarrow $ and up
diagonal step $\searrow $. Let $g(p,q)$ be the number of $HVD$-lattice paths
from $(0,0)$ to $(p,q)$.

We can analyze HVD-lattice path step by step. For each point on $x$-axis or $%
y$-axis, there is only one path from $(0,0)$ to them, it only $\rightarrow $
or $\downarrow $, namely $g(k,0)=g(0,r)=1$. For the point $(1,1)$, it is
either from the left point $(1,0)$ by $\rightarrow $ or from the above point 
$(0,1)$ by $\downarrow $, or from the point (0,0) by $\searrow $ to it, thus
the number $g(1,1)$ is the sum of $g(1,0),$ $\ g(0,1)$, and $g(0,0)$, that
is $g(1,1)=g(1,0)+g(0,1)+g(0,0)$. In general for the point $(k,r)$, the last
step of the $HVD$-lattice path is either from $(k,r-1)$ by $\rightarrow $ or
from $(k-1,r)$ by $\downarrow $, or from the $(k-1,r-1)$ by $\searrow $\ to
it, so there is 
\begin{equation}
g(k,r)=g(k,r-1)+g(k-1,r)+g(k-1,r-1).  \label{gg1}
\end{equation}%
this is as same as in paper \cite{Square}. In \cite{www2}, $g(k,r)$ is
called as \textit{Delanny number}. By tabulating the number of $HVD$-lattice
paths from $(0,0)$ to each point, we obtain a symmetric number table (see
table 3) which is as same as the Table 1 in paper \cite{Square} and
symmetric matrix $K$ as follows.

$%
\begin{tabular}{c|cccccccc}
\multicolumn{9}{l}{Table 3} \\ \hline
$p\backslash q$ & $0$ & $1$ & $2$ & $3$ & $4$ & $5$ & $6$ & $\cdots $ \\ 
\hline
$0$ & $1$ & $1$ & $1$ & $1$ & $1$ & $1$ & $1$ & $\cdots $ \\ 
$1$ & $1$ & $3$ & $5$ & $7$ & $9$ & $11$ & $13$ & $\cdots $ \\ \cline{6-7}
$2$ & $1$ & $5$ & $13$ & $25$ & \multicolumn{1}{|c}{$41$} & $61$ & 
\multicolumn{1}{|c}{$85$} & $\cdots $ \\ 
$3$ & $1$ & $7$ & $25$ & $63$ & \multicolumn{1}{|c}{$129$} & $231$ & 
\multicolumn{1}{|c}{$377$} & $\cdots $ \\ \cline{6-7}
$4$ & $1$ & $9$ & $41$ & $129$ & $321$ & $681$ & $1289$ & $\cdots $ \\ 
$5$ & $1$ & $11$ & $61$ & $231$ & $681$ & $1683$ & $3653$ & $\cdots $ \\ 
$6$ & $1$ & $13$ & $85$ & $377$ & $1289$ & $3653$ & $8989$ & $\cdots $ \\ 
$\vdots $ & $\vdots $ & $\vdots $ & $\vdots $ & $\vdots $ & $\vdots $ & $%
\vdots $ & $\vdots $ & 
\end{tabular}%
$

$K=\left[ 
\begin{tabular}{cccccccc}
$1$ & $1$ & $1$ & $1$ & $1$ & $1$ & $1$ & $\cdots $ \\ 
$1$ & $3$ & $5$ & $7$ & $9$ & $11$ & $13$ & $\cdots $ \\ 
$1$ & $5$ & $13$ & $25$ & $41$ & $61$ & $85$ & $\cdots $ \\ 
$1$ & $7$ & $25$ & $63$ & $129$ & $231$ & $377$ & $\cdots $ \\ 
$1$ & $9$ & $41$ & $129$ & $321$ & $681$ & $1289$ & $\cdots $ \\ 
$1$ & $11$ & $61$ & $231$ & $681$ & $1683$ & $3653$ & $\cdots $ \\ 
$1$ & $13$ & $85$ & $377$ & $1289$ & $3653$ & $8989$ & $\cdots $ \\ 
$\vdots $ & $\vdots $ & $\vdots $ & $\vdots $ & $\vdots $ & $\vdots $ & $%
\vdots $ & 
\end{tabular}%
\right] $

\begin{thm}
\cite{Richard} Let $r\leq \min \{p,q\}$. Then%
\begin{equation}
g(p,q)=\overset{\min \{p,q\}}{\underset{r=0}{\sum }}\left( 
\begin{array}{ccc}
& p+q-r &  \\ 
p-r & q-r & r%
\end{array}%
\right) =\underset{r=0}{\overset{\min \{p,q\}}{\sum }}\frac{(p+q-r)!}{%
(p-r)!(q-r)!r!}  \label{r6}
\end{equation}

where $\left( 
\begin{array}{ccc}
& p+q-r &  \\ 
p-r & q-r & r%
\end{array}%
\right) $ is the multinomial number \cite{Richard}.
\end{thm}

\begin{thm}
\cite{Square}%
\begin{equation}
g(p,q)=\underset{\alpha =0}{\overset{p}{\sum }}\binom{p}{\alpha }\binom{%
q+\alpha }{p}  \label{r66}
\end{equation}
\end{thm}

\begin{thm}
\cite{Square}%
\begin{equation}
g(p,q)=\underset{\alpha =0}{\overset{\min \{p,q\}}{\sum }}\binom{p}{\alpha }%
\binom{q}{\alpha }2^{\alpha }  \label{r666}
\end{equation}
\end{thm}

We can rewrite (\ref{r666}) to matrix format, there is the following
proposition.

\begin{prop}
The first $n$ rows and $n$ columns $K_{n}$ is the product of Pascal matrix $%
P_{n}$, diagonal matrix $D_{n}$ and its transposition of Pascal matrix $%
P_{n}^{T}$.%
\begin{equation}
K_{n}=P_{n}D_{n}P_{n}^{T}  \label{c1}
\end{equation}%
where $D_{n}=diag(1,2,2^{2},2^{3},\cdots ,2^{n-1}).$
\end{prop}

\section{Subdiagonal HVD-lattice path and Schr\H{o}der number}

We now consider subdiagonal $HVD$-lattice path from $(0,0)$ to $(p,q)$ that
is restricted to lie on or below the line $y=x$ in the coordinate plane
whenever $p\geq q$. Let $R(p,q)$ be the number of the subdiagonal $HVD$%
-lattice path from $(0,0)$ to $(p,q)$.

When we analyze subdiagonal $HVD$-lattice path step by step, we obtain a
triangle (see table 4)

$%
\begin{tabular}{c|ccccccccc}
\multicolumn{10}{l}{Table 4} \\ \hline
$p\backslash q$ & $0$ & $1$ & $2$ & $3$ & $4$ & $5$ & $6$ & $7$ & $\cdots $
\\ \hline
$0$ & $1$ &  &  &  &  &  &  &  &  \\ 
$1$ & $1$ & $2$ &  &  &  &  &  &  &  \\ 
$2$ & $1$ & $4$ & $6$ &  &  &  &  &  &  \\ 
$3$ & $1$ & $6$ & $16$ & $22$ &  &  &  &  &  \\ \cline{5-6}\cline{5-6}
$4$ & $1$ & $8$ & $30$ & \multicolumn{1}{|c}{$68$} & $90$ & 
\multicolumn{1}{|c}{} &  &  &  \\ 
$5$ & $1$ & $10$ & $48$ & \multicolumn{1}{|c}{$146$} & $304$ & 
\multicolumn{1}{|c}{$394$} &  &  &  \\ \cline{5-6}\cline{5-6}
$6$ & $1$ & $12$ & $70$ & $264$ & $714$ & $1412$ & $1806$ &  &  \\ 
$7$ & $1$ & $14$ & $96$ & $430$ & $1408$ & $3534$ & $6752$ & $8558$ &  \\ 
$\vdots $ & $\vdots $ & $\vdots $ & $\vdots $ & $\vdots $ & $\vdots $ & $%
\vdots $ & $\vdots $ & $\vdots $ & 
\end{tabular}%
$

Compare the lower triangle of Table 3 and Table 4, we find that the element
in Table 4 $R(p,q)$ equal the element in Table 3 $g(p,q)$ subtract the
element $g(p+1,q-1)$, that is $R(p,q)=g(p,q)-g(p+1,q-1)$. This is derived
from the symmetric of the Table 3, $g(p+1,q-1)$ is the number of subdiagonal 
$HVD$-lattice paths above the line $y=x$, it is sure that it must be
subtracted. By Theorem 6, 7, 8, we can get

\begin{thm}
Let $p$ and $q$ be positive integers with $q\geq q$, and let $r$ be a
nonnegative integer with $r\leq q$. Then%
\begin{eqnarray}
R(p,q) &=&g(p,q)-g(p+1,q-1)  \label{h11} \\
&=&\overset{\min \{p,q\}}{\underset{r=0}{\sum }}\left( 
\begin{array}{ccc}
& p+q-r &  \\ 
p-r & q-r & r%
\end{array}%
\right) -\overset{\min \{p+1,q-1\}}{\underset{r=0}{\sum }}\left( 
\begin{array}{ccc}
& p+q-r &  \\ 
p+1-r & q-1-r & r%
\end{array}%
\right) ,  \notag
\end{eqnarray}%
\begin{eqnarray}
R(p,q) &=&g(p,q)-g(p+1,q-1)  \label{h2} \\
&=&\underset{r=0}{\overset{p}{\sum }}\binom{p}{r}\binom{q+r}{p}-\underset{r=0%
}{\overset{p+1}{\sum }}\binom{p+1}{r}\binom{q-1+r}{p+1},  \notag \\
&&  \notag
\end{eqnarray}%
\begin{eqnarray}
R(p,q) &=&g(p,q)-g(p+1,q-1)  \label{h3} \\
&=&\underset{r=0}{\overset{\min \{p,q\}}{\sum }}\binom{p}{r}\binom{q}{r}%
2^{r}-\underset{r=0}{\overset{\min \{p+1,q-1\}}{\sum }}\binom{p+1}{r}\binom{%
q-1}{r}2^{r},  \notag
\end{eqnarray}%
\begin{equation}
R(p,q)=\underset{r=0}{\overset{q}{\sum }}\frac{p-q+1}{p-r+1}\left( 
\begin{array}{ccc}
& p+q-r &  \\ 
p-r & q-r & r%
\end{array}%
\right) .\text{ \ \ \ \ \ \ \ \ \ \ \ \ \ \ \ \ \ \ \ \ \ \ \ \ \ }
\label{h4}
\end{equation}%
where (\ref{h4}) is in \cite{Richard}.
\end{thm}

We now suppose that $p=q=n$. The subdiagonal HVD-lattice paths from $(0,0)$
to $(n,n)$ are called Schr\H{o}der paths [1, 16, 17]. The \textit{large Schr%
\H{o}der number} $S_{n}$ is the number of Schr\H{o}der paths from $(0,0)$ to 
$(n,n)$. Thus, we have%
\begin{equation}
S_{n}=R(n,n)=\underset{r=0}{\overset{n}{\sum }}\frac{1}{n-r+1}\left( 
\begin{array}{ccc}
& 2n-r &  \\ 
n-r & n-r & r%
\end{array}%
\right) ,  \label{gf1}
\end{equation}%
\begin{equation}
S_{n}=R(n,n)=2^{n}+\underset{r=0}{\overset{n-1}{\sum }}[\frac{r2^{r}}{n-r+1}%
\binom{n}{r}\binom{n-1}{r}].  \label{gf2}
\end{equation}

By the formula (\ref{h11}-\ref{h4}), we obtain following proposition.

\begin{prop}
The numbers $R(p,q)$ are given by the recurrence relation 
\begin{equation}
R(p,q)=R(p-1,q-1)+R(p-1,q)+R(p,q-1),  \label{s1}
\end{equation}%
\begin{equation}
S_{n}=R(n,n)=R(n-1,n-1)+R(n,n-1),  \label{s11}
\end{equation}%
\begin{equation}
R(p,q)=2[R(p-1,0)+R(p-1,1)+\cdots +R(p-1,q-1)]+R(p-1,q),\text{if }p\neq q,
\label{s2}
\end{equation}%
\begin{equation}
R(n,n)=2[R(n-1,0)+R(n-1,1)+\cdots +R(n-1,n-1)],  \label{s31}
\end{equation}%
\begin{equation}
S_{n}=S_{n-1}+S_{0}S_{n-1}+S_{1}S_{n-2}+S_{2}S_{n-3}+\cdots +S_{n-1}S_{0}
\label{s311}
\end{equation}%
\begin{equation}
(n+2)S_{n+2}=3(2n+1)S_{n+1}-(n-1)S_{n}.  \label{s32}
\end{equation}%
\bigskip
\end{prop}

\end{document}